\documentclass[aap,MSNbibl,dvips]{arximspdf}

\doi{10.1214/09-AAP658}
\volume{20}
\issue{3}
\pubyear{2010}
\firstpage{890}
\lastpage{906}

\makeatletter
\newtheorem{theorem}{Theorem}

\newtheorem{corollary}[theorem]{Corollary}

\newtheorem{lemma}[theorem]{Lemma}

\newtheorem{proposition}[theorem]{Proposition}

\newproclaim{definition}[theorem]{Definition}
\newproclaim{assumption}{Assumption}
\makeatother

\begin{document}
\begin{frontmatter}

\title{Products of random matrices: Dimension and growth in norm}
\runtitle{Products of random matrices}

\begin{aug}
\author[A]{\fnms{Vladislav}  \snm{Kargin}\corref{}\ead[label=e1]{kargin@stanford.edu}}
\runauthor{V. Kargin}
\affiliation{Stanford University}
\address[A]{Department of Mathematics\\
Stanford University\\
282 Mosher Way\\ Palo Alto, California 94304\\
USA\\
\printead{e1}} 
\end{aug}

\received{\smonth{1} \syear{2009}}

\begin{abstract}
Suppose that $X_1,\ldots ,X_n,\ldots $ are i.i.d. rotationally invariant $N$-by-$N$ matrices. Let
$\Pi _n=X_n \cdots X_1$. It is known that $n^{-1} \log |\Pi_n| $ converges to a nonrandom limit. We prove that under
certain additional assumptions on matrices $X_i$ the speed of convergence
to this limit does not decrease when the size of matrices, $N$, grows.
\end{abstract}



\end{frontmatter}

\section{Introduction}

Let $X_{i}$ be a sequence of independent $N\times N$ random matrices and $%
\Pi _{n}=X_{n}\cdots X_{1}.$ In a celebrated paper \cite%
{furstenbergkesten60}, Furstenberg and Kesten proved that $n^{-1}\log
\| \Pi _{n}\| $ converges provided that $E\log ^{+}( \|
X_{i}\| ) <\infty .$ Later, Oseledec in \cite{oseledec68} proved
convergence for other singular values of $\Pi _{n},$ and Cohen and Newman in %
\cite{cohennewman84} studied the behavior of the limit in the situation
when $N$ approaches infinity. This paper investigates the question of how
the speed of convergence depends on the dimension of matrices $N$.

Consider a dynamical system (a gas, an economy, an ecosystem, etc.). Its
evolution can be described by a mapping $\psi _{i}\rightarrow X_{i}(
\psi _{i}) ,$ where $\psi _{i}$ is a vector that describes the state
of the system at time $i.$ We can often model the mapping as a
multiplication by a random matrix $X_{i}.$ Stability and other long-run
properties of the system depend on the growth in the norm of the product $%
\Pi _{n}=X_{n}\cdots X,$ which we can measure by calculating the quantity $%
n^{-1}\log ( \| \Pi _{n}\| ) .$

The sub-multiplicativity property of the norm ($\| X_{2}X_{1}\|
\leq \| X_{2}\| \| X_{1}\| $) ensures that $n^{-1}\log
( \| \Pi _{n}\| ) $ converges to $E\log \|
X_{1}u\| ,$ where $u$ is an arbitrary vector. Intuitively, this means
that it is not important what was the starting vector of the system. After
some time, all products grow at the same rate independently of the initial
state.

It is of interest to investigate whether this erasure of memory about the
initial state occurs slower in more complex systems, that is, in systems,
which are described by matrices of larger size.

Of course, when we compare long-run properties of systems, we should only
look at the systems that are comparable in the short run, that is, the
system that have comparable one-step behavior. Roughly, the difference
between one-step growth of a specially-chosen and a random vector can be
measured by the ratio of $\| X_{1}\| ^{2}$ to $N^{-1}\mathrm{tr}%
( X_{1}^{\ast }X_{1}) ,$ where $N$ is the dimension of the matrix
$X_{i}. $ Indeed, $\| X_{1}\| ^{2}$ is the square of the maximal
possible increase in the length of the state vector, and $N^{-1}\operatorname{tr}%
( X_{1}^{\ast }X_{1}) $ is the average of the squared singular
values of $X_{1},$ hence it can be considered as a measure of the increase
in the length of a random state vector.

Hence, if we want systems to be comparable in the short run, then we should
restrict this ratio by a constant that does not depend on the dimension of
the system. (Otherwise, some directions may become more and more unusual as
the dimension of the system grows.)  We will call this property uniform
boundedness of singular values.

We also want to look at sufficiently symmetric systems, that is, systems
without preferential directions. We codify this by requiring that matrices $%
X_{i}$ are rotationally invariant, that is, the distribution of matrix
elements does not depend on the choice of basis.

The main result of this paper is that under these assumptions the speed
with which the memory of the initial state is erased does not decrease as
the dimension of the system grows.

Intuitively, the asymptotic behavior of $n^{-1}\log \| \Pi _{n}\|
$ depends on three factors. First of all, for a fixed vector $v,$%
\[
n^{-1}\log \| \Pi _{n}v\| =n^{-1}\sum_{i=1}^{n}\log \|
X_{i}v_{i}\|
\]%
for a certain sequence of vectors $v_{i}$ and this averaging is likely to
concentrate the distribution of $n^{-1}\log \| \Pi _{n}v\| .$
This factor does not depend on the dimension $N$. On the other hand, we are
interested in the convergence of the supremum of $n^{-1}\log \| \Pi
_{n}v\| $ over all $v\in S^{N},$ and to ensure the convergence of this
supremum we have to make sure that variables $n^{-1}\log \| \Pi
_{n}v\| $ are all close to the limit $E\log \| X_{1}u\| $
for a sufficiently dense set of vectors $v.$ The number of elements in such
a set is likely to grow exponentially in $N,$ and this might make the
convergence of $n^{-1}\log \| \Pi _{n}\| $ slower for large $N.$

The third factor appears because for every fixed vector $v,$ the norm $%
\Vert X_{i}v\Vert $ becomes concentrated around some particular
value as $N\rightarrow \infty .$ This factor is likely to speed up the
convergence of $n^{-1}\log \Vert \Pi _{n}v\Vert $ and therefore
of $n^{-1}\log \Vert \Pi _{n}\Vert .$

We will show in this paper that the third factor dominates and the speed of
convergence of $n^{-1}\log \| \Pi _{n}\| $ is not slowed down by
the growth in the dimension~$N$.

Previously, the speed of convergence in the Furstenberg--Kesten theorem was
investigated in \cite{lepage82,tutubalin65} and \cite%
{guivarchraugi85}. They proved a central limit theorem for $n^{-1/2}\log
\| \Pi _{n}\| $ and studied large deviations of $n^{-1}\log
\| \Pi _{n}\| $ for a large class of random matrices. However,
the results in these papers do not provide effectively computable bounds on
the rate of convergence in limit theorems, and, as a consequence, do not
help us to investigate how the speed of convergence changes as the dimension
of matrices grows. One of the contributions of this paper is deriving more
explicit bounds on the speed of convergence in limit theorems.

Let us describe the problem in a more formal fashion. Consider independent
identically-distributed $N$-by-$N$ matrices $X_{i}^{( N) }.$ We
are interested in the behavior of the norm of the product $\Pi
_{n}=X_{n}^{( N) }\cdots X_{1}^{( N) },$ and we will
make the following assumptions about matrices $X_{i}^{( N) }$.
First of all, we assume that random matrices $X_{i}^{( N) }$ are
rotationally invariant; that is, the distribution of their entries does not
depend on the choice of coordinates. Formally, we use the following
definition.

\begin{definition}
A random matrix $X$ is rotationally invariant if for every integer $k\geq 1,$
for every collection of vectors $\{ v_{i},w_{i}\} $, $i=1,\ldots,k,$
and for every orthogonal matrix $U,$ the joint distributions of random
vectors $\{ \langle w_{i},Xv_{i}\rangle \} _{i=1}^{k}$
and $\{ \langle Uw_{i},XUv_{i}\rangle \} _{i=1}^{k}$
are the same.
\end{definition}

\renewcommand{\theassumption}{\Alph{assumption}}
\setcounter{assumption}{0}
\begin{assumption}[(``Rotational invariance'')]\label{AA} Matrices $X_{i}^{(
N) }$ are rotationally invariant.
\end{assumption}

We also impose an assumption needed for the validity of the
Furstenberg--Kesten theorem.

\begin{assumption}[(``Furstenberg--Kesten'')]\label{AB} For all $N,$ $E\log
^{+}\| X_{i}^{( N) }\| $ exists.
\end{assumption}

Second, we restrict our study to two important cases. The first one is the
case of (real) Gaussian matrices $X_{i}^{( N) }$, that is,
independent random $N$-by-$N$ matrices with independent entries distributed
according to the Gaussian distribution with zero mean and variance $\sigma
^{2}/N,$ that is, as $\mathcal{N}( 0,\sigma ^{2}/N) .$

The second case is that of independent rotationally invariant $N$-by-$N$
matrices $X_{i}^{( N) }$ that satisfy the following assumptions.
Let $s_{k}^{( i,N) }$ be the eigenvalues of $X_{i}^{(
N) \ast }X_{i}^{( N) }$ (i.e., squared singular values of $%
X_{i}^{( N) }$), and let
\[
\overline{s}^{( i,N) }=\frac{1}{N}\sum_{k-1}^{N}s_{k}^{(
i,N) }=\frac{1}{N}\operatorname{tr}\bigl( X_{i}^{( N) \ast
}X_{i}^{( N) }\bigr) .
\]%
(We will sometimes omit superscripts to lighten the notation.)

\begin{assumption}[(``Uniformly bounded singular values'')]\label{AC} With probability
1, $\max_{k}s_{k}^{( i,N) }\leq b\overline{s}^{( i,N)
},$ where the constant $b$ does not depend on $N.$
\end{assumption}

In other form, Assumption~\ref{AC} says that
\[
\bigl\| X_{i}^{( N) }\bigr\| ^{2}\leq b\frac{1}{N}\operatorname{tr}%
\bigl( X_{i}^{( N) \ast }X_{i}^{( N) }\bigr)
\]
with probability 1.

\begin{assumption}[(``Comparability across $N$'')]\label{AD} \hspace*{-2pt}$\operatorname{var}[ \log
\overline{s}^{( i,N) }] $ exists and bound\-ed by a constant
which does not depend on $N.$
\end{assumption}

One example of a matrix family that satisfies these assumptions is Hermitian
matrices $X_{i}^{( N) }$ which are generated in the following
way. Sample $N$ independent values from a distribution supported on $[
\alpha ,\beta ] ,$ where $\beta >\alpha >0,$ and construct a diagonal
matrix $D^{( N) }$ by putting these values on the main diagonal.
Then take a Haar-distributed random orthogonal matrix $U_{i}^{(
N) }$ and define $X_{i}^{( N) }$ as $D^{( N)
}U_{i}^{( N) }.$ A sequence of these matrices (with independent $%
U_{i}^{( N) }$) will satisfy all the assumptions.

The main result is as follows.

\begin{theorem}
\label{theorem_uniform_convergence_norm_products}Let $X_{i}^{( N)
}$ be independent, identically distributed $N\times N$ matrices, which
satisfy Assumptions \ref{AA} and \ref{AB} and which are either Gaussian with
independent entries $\mathcal{N}( 0,\sigma ^{2}/N) ,$ or satisfy
Assumptions \ref{AC} and \ref{AD}. Let $\Pi _{n}=X_{n}^{( N) }\cdots
X_{1}^{( N) }$ and let $v$ be an arbitrary unit vector. Then $%
n^{-1}\log \| \Pi _{n}\| $\vspace*{-1pt} converges in probability to $E\log
\| X_{1}^{( N) }v\| $ and the convergence is uniform
in $N.$ That is, for each $\delta >0,$ there exists an $n_{0}( \delta
) $ such that for all $n\geq n_{0}$ and all $N\geq 1,$
\begin{equation}\label{inequality_theorem}
\Pr \bigl\{ \bigl| n^{-1}\log \| \Pi _{n}\| -E\log \bigl\|
X_{1}^{( N) }v\bigr\| \bigr| \geq \delta \bigr\} \leq \delta .
\end{equation}
\end{theorem}

The assumptions of the theorem are sufficient but not necessary. The
assumption that $s_{k}\leq b\overline{s}$ is used in the proof of
Proposition \ref{Proposition_LD_inequality_unbounded_rv} below, where it is
used to estimate the probability of large deviations of $\log \|
X_{i}^{( N) }v\| $ and to show that the rate in the
corresponding exponential inequality is proportional to $N.$ It is likely
that this assumption can be somewhat relaxed by requiring instead that $\Pr
\{ s_{k}/\overline{s}>b+u\} \leq ce^{-c^{\prime }Nu}$.

One particular implication of the assumption $s_{k}\leq b\overline{s}$ is
that the bound on singular values does not depend on the dimension of the
matrix. In order to understand this assumption better, consider the
following example. Let
\[
X_{i}^{( N) }=\sqrt{N}| y_{i}\rangle \langle
x_{i}| ,
\]%
where $\langle x_{i}| $ is a Haar-distributed row $N$-dimensional
vector, and $| y_{i}\rangle $ is a Haar-distributed column $N$%
-dimensional vector. (Vectors $\langle x_{i}| $ and $|
y_{i}\rangle $ are assumed to be independent.) Then the squared
singular values of $X_{i}$ are all zero except one, which equals $N.$ Hence, $%
\overline{s}^{( i,N) }=1$ and $\log \overline{s}^{(
i,N) }=0.$ We can conclude that Assumptions \ref{AA}, \ref{AB} and \ref{AD} are
satisfied, and Assumption \ref{AC} is not satisfied.

Next, consider
\[
\xi _{i}:=\bigl\| X_{i}^{( N) }v \bigr\| ^{2}=N\langle
x_{i}|v\rangle ^{2},
\]
where $v$ is an arbitrary vector. It is easy to see that $\xi _{i}$ is
distributed as the first coordinate of a Haar-distributed vector $u.$ In
other words $\xi _{i}$ is distributed as
\[
\frac{Y_{1}^{2}}{( Y_{1}^{2}+\cdots + Y_{N}^{2}) /N},
\]%
where $Y_{i}$ are independent standard Gaussian variables. Then it is clear
that
\[
\lim_{N\rightarrow \infty }E\log \| X_{1}v\| ^{2}=E\log (
Y_{1}) ^{2}\in ( -\infty ,0) .
\]

Next, let us compute $n^{-1}\log \| \Pi _{n}\| ^{2}.$ Note that
\[
\Pi _{n}=N^{n/2}| y_{n}\rangle \langle x_{n}| \cdots
| y_{2}\rangle \langle x_{2}|y_{1}\rangle \langle
x_{1}|
\]%
and
\begin{eqnarray*}
\Pi _{n}^{\ast }\Pi _{n} &=&N^{n}| x_{1}\rangle \langle
x_{n}|y_{n-1}\rangle ^{2}\cdots \langle x_{2}|y_{1}\rangle
^{2}\langle x_{1}|
\\
&=&N| x_{1}\rangle N\langle x_{n}|y_{n-1}\rangle
^{2}\cdots N\langle x_{2}|y_{1}\rangle ^{2}\langle
x_{1}| .
\end{eqnarray*}
Hence,
\[
n^{-1}\log \| \Pi _{n}\| ^{2}=\frac{\log N}{n}+n^{-1}%
\sum_{i=1}^{n-1}\log \xi _{i},
\]%
where $\xi _{i}$ are independent and distributed as $\| X_{i}^{(
N) }v\| ^{2}$ above. Hence, $\xi _{i}$ converges in distribution
to $Y_{1}^{2}$ as $N\rightarrow \infty .$ It is clear that
\[
n^{-1}\sum_{i=1}^{n-1}\log \xi _{i}\rightarrow E\log \| X_{1}v\|
^{2}
\]%
in probability as $n\rightarrow \infty .$ Therefore, for large $N,$
\[
n^{-1}\log \| \Pi _{n}\| ^{2}-E\log \| X_{1}v\|
^{2}\sim \frac{\log N}{n}.
\]%
This bias term cannot be made small uniformly in $N$ by an increase in $n.$
This means that the claim of Theorem \ref%
{theorem_uniform_convergence_norm_products} fails in this case.

Later, in Section \ref{section_necessity}, we will prove a necessary
condition for the uniform convergence by using the basic idea of this
example.

In order to understand the role of the rotational invariance assumption,
consider the following example.

Let $X_{i}$ be independent, identically distributed, diagonal matrices. The
diagonal elements of a matrix $X_{i}$ are independent Bernoulli variables
that take values $a$ and $b$. That is, a diagonal element takes the value $%
b>0$ with probability $p$ and the value $a>0$ with probability $q=1-p.$
Assume that $b>a.$

It is easy to see that the norm of $\Pi _{n}=X_{1}\cdots X_{n}$ is given by
the following expression:%
\[
\| \Pi _{n}\| =\max \{ a^{\alpha _{1}}b^{\beta _{1}},\ldots
,a^{\alpha _{N}}b^{\beta _{N}}\} ,
\]%
where $\alpha _{i}+\beta _{i}=n,$ and $\beta _{i}$ are independent random
variables with the binomial distribution $B( p,n).$

Taking the logarithm and dividing by $n,$ we get
\[
\frac{1}{n}\log \Vert \Pi _{n}\Vert =\log a+\log (
b/a) \max \{ \widetilde{\beta }_{1},\ldots ,\widetilde{\beta }%
_{N}\},
\]%
where $\widetilde{\beta }_{i}=\beta _{i}/n.$ Note that as $n$ grows, each $%
\widetilde{\beta }_{i}$ approaches the Gaussian distribution $\mathcal{N}%
( p,pq/n) .$

If $N$ is fixed, then $\lim n^{-1}\log \| \Pi _{n}\| =\log
a+p\log ( b/a) .$ However, if $N$ grows simultaneously with $n,$
then the limit of $n^{-1}\log \| \Pi _{n}\| $ may be
nonexistent, or may depend on the speed of growth in $N$ relative to the
speed of growth in $n.$ Hence, the conclusion of Theorem \ref%
{theorem_uniform_convergence_norm_products} is invalid in this case.

It is an interesting problem whether the assumption of rotational invariance
can be relaxed so that the result in Theorem \ref%
{theorem_uniform_convergence_norm_products} holds for a larger class of
matrices, for example, for matrices with i.i.d. non-Gaussian entries (i.e.,
Wigner matrices). However, this problem appears to be hard since at this
moment very little is known about effective bounds on the rate of
convergence in the Furstenberg--Kesten theorem.

Let me now explain two results which will be used as tools in the proof of
Theorem \ref{theorem_uniform_convergence_norm_products}. The proofs of these
results will be given in later sections.

Our main tool is the following proposition.

\begin{proposition}
\label{Proposition_LD_inequality_unbounded_rv} \textup{(i)} Suppose that all $X_{i}$
are Gaussian with independent entries $\mathcal{N}( 0,\sigma
^{2}/N) $. Then for all sufficiently small $t$, all $N\geq N_{1}(
t) $ and all $n\geq 1,$
\[
\Pr \biggl\{ \biggl| \frac{1}{n}\log \| \Pi _{n}v\| -\log \sigma
\biggr| >t\biggr\} \leq 2\exp \biggl( -\frac{1}{8}Nnt^{2}\biggr) .
\]
\textup{(ii)} Suppose that i.i.d. $N$-by-$N$ matrices $X_{i}$ are rotationally
invariant and satisfy Assumption~\ref{AC} with constant $b.$ Let
\[
\overline{s}^{( i,N) }=\frac{1}{N}\sum_{k-1}^{N}s_{k}^{(
i,N) }=\frac{1}{N}\mathrm{Tr}( X_{i}^{\ast }X_{i}) .
\]%
Then for all $t\in ( 0,1/4) $, all $N\geq N_{1}( t) $
and all $n\geq 1,$
\[
\Pr \Biggl\{ \Biggl| \frac{1}{n}\log \| \Pi _{n}v\| -\frac{1}{n}%
\sum_{i=1}^{n}\log \overline{s}^{( i,N) }\Biggr| >t\Biggr\} \leq
2\exp \biggl( -\frac{1}{32b^{2}}Nnt^{2}\biggr) .
\]
\end{proposition}

In its essence, Proposition \ref{Proposition_LD_inequality_unbounded_rv} is
a large deviation result which quantifies the speed of convergence of $%
n^{-1}\log \Vert \Pi _{n}v\Vert ^{2}$ for a fixed vector $v.$ Its
main point is that the rate in this large deviation estimate is proportional
to the dimension $N$. The proof of this proposition will be given in Section %
\ref{section_large_deviation_bound}.

The other tool is as follows. Let a set of points on the unit sphere in $%
\mathbb{R}^{N}$ be called an $\varepsilon  $-net if the sphere is covered
by spherical caps with centers at these points and angular radius $%
\varepsilon .$

\begin{proposition}
\label{Proposition_norm_and_sup}Let $A$ be an arbitrary $N$-by-$N$ matrix.
Suppose that the endpoints of vectors $v_{i}$ form an $\varepsilon $-net of
the unit sphere in $\mathbb{R}^{N}$. Then for all sufficiently small $%
\varepsilon $,
\[
\log \| A\| \leq \max_{i}\log \| Av_{i}\|
+2\varepsilon .
\]
\end{proposition}

This proposition allows us to control the matrix norm $\| \Pi
_{n}\| $ by the norms of vectors $\| \Pi _{n}v_{i}\| ,$
where $v_{i}$ runs through a finite set of values.

\begin{pf*}{Proof of Proposition \protect\ref{Proposition_norm_and_sup}} Let $v_{i}$ be a vector in the net which is closest to a
unit vector $v.$ Then%
\begin{eqnarray*}
\| Av\| &\leq &\| Av_{i}\| +\| A(
v-v_{i}) \|
\\
&\leq &\| Av_{i}\| +\varepsilon \| A\| .
\end{eqnarray*}%
Taking the supremum over $v,$ we obtain that
\[
( 1-\varepsilon ) \| A\| \leq \max_{i}\|
Av_{i}\| .
\]%
Hence,
\[
\log \| A\| \leq \max_{i}\log \| Av_{i}\| -\log (
1-\varepsilon ) ,
\]%
and the claim of the proposition follows.
\end{pf*}

This proposition is useful in conjunction with the following result about
the size of sphere coverings. By Lemma 2.6 on page 7 of \cite%
{milmanschechtman01}, for $\varepsilon $ smaller than a certain constant,
there exists an $\varepsilon $-net with cardinality $M\leq \exp ( N\log
( 3/\varepsilon ) ) .$

Now let us prove Theorem \ref{theorem_uniform_convergence_norm_products} by
using Propositions \ref{Proposition_LD_inequality_unbounded_rv} and \ref%
{Proposition_norm_and_sup}.

\begin{pf*}{Proof of Theorem \protect\ref{theorem_uniform_convergence_norm_products}}
We focus on the case when Assumptions \ref{AC} and \ref{AD} hold. The proof for the
case of Gaussian matrices goes along a similar route and it is simpler.

First of all, note that is enough to prove that (\ref{inequality_theorem})
holds for all sufficiently large $N,$ that is, for all $N\geq N_{0}(
\delta ) .$ Indeed, for each $N\leq N_{0}$ we can apply results in %
\cite{furstenberg63} and find that inequality (\ref{inequality_theorem})
holds if $n\geq n( \delta ,N) .$ Hence, inequality (\ref%
{inequality_theorem}) holds for all $N\leq N_{0},$ provided that
\[
n\geq n_{0}( \delta ) =\max_{N\leq N_{0}( \delta )
}\{ n( \delta ,N) \} .
\]

We will choose an appropriate $N_{0}( \delta ) $ later.

We are going to prove that for all sufficiently large $N$ and $n$ [i.e., $%
N\geq N_{2}( \delta ) $ and all $n\geq n_{2}( \delta )
$], it is true that
\begin{equation}\label{inequality_proof_3}
\Pr \Biggl\{ \Biggl| \frac{1}{n}\log \| \Pi _{n}\| ^{2}-\frac{1}{n}%
\sum_{i=1}^{n}\log \overline{s}^{( i,N) }\Biggr| >\frac{\delta }{%
10}\Biggr\} <\frac{\delta }{10}.
\end{equation}
Let vectors $v_{j},$ $j=1,\ldots ,M,$ form a $( \delta /100) $%
-net on the unit sphere. Then, by using Propositions \ref%
{Proposition_norm_and_sup} and \ref{Proposition_LD_inequality_unbounded_rv},
the union bound and the estimate on the number of elements in the net we
obtain
\begin{eqnarray*}
&&\Pr \Biggl\{ \Biggl| \frac{1}{n}\log \| \Pi _{n}\| ^{2}-\frac{1}{n%
}\sum_{i=1}^{n}\log \overline{s}^{( i,N) }\Biggr| >\frac{\delta }{%
10}\Biggr\}
 \\
&&\qquad \leq \Pr \Biggl\{ \max_{v_{i}}\Biggl| \frac{1}{n}\log \| \Pi
_{n}v_{i}\| ^{2}-\frac{1}{n}\sum_{i=1}^{n}\log \overline{s}^{(
i,N) }\Biggr| >\frac{\delta }{100}\Biggr\}
 \\
&&  \qquad\leq 2\exp \biggl\{ \biggl( \log \biggl( \frac{300}{\delta }\biggr) -cn\biggl(
\frac{\delta }{100}\biggr) ^{2}\biggr) N\biggr\} ,
\end{eqnarray*}
where $c$ is a certain constant. Clearly, we can choose $n_{2}( \delta
) $ in such a way that for all $n\geq n_{2}( \delta ) ,$ it
is true that
\[
\log \biggl( \frac{300}{\delta }\biggr) -cn\biggl( \frac{\delta }{100}\biggr)
^{2}<\alpha <0
\]%
for some $\alpha ,$ and then choose $N_{2}( \delta ) ,$ such that
for all $N>N_{2}( \delta ) $ it is true that%
\[
2\exp \{ \alpha N\} <\frac{\delta }{10}.
\]
This choice of $n_{2}( \delta ) $ and $N_{2}( \delta )
$ is sufficient to ensure that (\ref{inequality_proof_3}) holds.

Next, let $d_{N}=E\log \overline{s}^{( i,N) }.$ Since variance of
$\log \overline{s}_{i}^{( N) }$ is bounded above by a finite
constant which does not depend on $N$ (Assumption \ref{AD}), therefore we can
find such $n_{3}( \delta ) $ that for all $n\geq n_{3}(
\delta ) ,$ it is true that
\begin{equation}\label{inequality_proof_2}
\Pr \Biggl\{ \Biggl| \frac{1}{n}\sum_{i=1}^{n}\log \overline{s}_{i}^{(
N) }-d_{N}\Biggr| >\frac{\delta }{100}\Biggr\} <\frac{\delta }{100}
\end{equation}
for all $N.$

It follows that for all $n\geq n_{4}( \delta ) $ and $N\geq
N_{2}( \delta ) ,$ it is true that
\begin{equation}\label{inequality_proof_4}
\Pr \biggl\{ \biggl| \frac{1}{n}\log \| \Pi _{n}\|
^{2}-d_{N}\biggr| >\frac{\delta }{5}\biggr\} <\frac{\delta }{5}.
\end{equation}
Note that by the Furstenberg--Kesten theorem,
\begin{equation}\label{inequality_proof_5}
\Pr \biggl\{ \biggl| \frac{1}{n}\log \| \Pi _{n}\| ^{2}-E\log
\bigl\| X_{i}^{( N) }u\bigr\| ^{2}\biggr| >\frac{\delta }{5}%
\biggr\} <\frac{\delta }{5}
\end{equation}
for all $n\geq n_{5}( \delta ,N) .$ This implies that for all $%
N\geq N_{2}( \delta ) ,$ there exists such $n,$ that both
inequalities (\ref{inequality_proof_4}) and (\ref{inequality_proof_5}) hold.
This implies that for all such $N,$ and for all $\delta <1,$ the following
inequality holds:
\begin{equation}\label{inequality_proof_6}
\biggl| d_{N}-E\log \bigl\| X_{i}^{( N) }u \bigr\| ^{2}\biggr| <%
\frac{2\delta }{5}.
\end{equation}
Otherwise, the sum of the events in (\ref{inequality_proof_4}) and (\ref%
{inequality_proof_5}) would cover all probability space and hence the sum of
probabilities in (\ref{inequality_proof_4}) and (\ref{inequality_proof_5}), $%
2\delta /5,$ would have to be greater then 1. This contradicts to the
assumption that $\delta <1$.

Inequalities (\ref{inequality_proof_3}), (\ref{inequality_proof_2}) and (\ref%
{inequality_proof_6}) imply that
\[
\Pr \biggl\{ \biggl| \frac{1}{n}\log \| \Pi _{n}\| ^{2}-E\log
\| X_{i}u\| ^{2}\biggr| >\delta \biggr\} <\delta
\]%
for all $n>n_{0}( \delta ) $ and $N>N_{0}( \delta ) ,$
where $n_{0}$ and $N_{0}$ are sufficiently large functions of $\delta $.
\end{pf*}

It remains to complete the proof by proving Proposition \ref%
{Proposition_LD_inequality_unbounded_rv}. We will do this in the next
section.

The rest of the paper consists of Section \ref{section_large_deviation_bound}%
, which is devoted to the proof of Proposition \ref%
{Proposition_LD_inequality_unbounded_rv}, Section \ref{section_necessity},
which gives a necessary condition for uniform convergence in
Furstenberg--Kesten theorem, and Section \ref{section_conclusion}, which
concludes.

\section{A large deviation bound for the dilation of a fixed vector} \label{section_large_deviation_bound}

Everywhere in this section, we assume that random matrices $X_{i}$ are
independent, identically distributed and rotationally invariant, and that $%
\Pi _{i}=X_{i}X_{i-1}\cdots X_{1}.$ Let us consider the following random
variables:
\[
y_{i}=\log \biggl( \frac{\| X_{i}\Pi _{i-1}v\| }{\| \Pi
_{i-1}v\| }\biggr) .
\]%
It is known (e.g., \cite{cohennewman84}) that the random variables $y_{i}$
are independent and identically distributed. Their distribution coincides
with the distribution of $\log ( \| X_{1}v\| ) ,$
where $v$ is an arbitrary unit vector.

\subsection{Gaussian matrices}

In this section, we consider an important case when each matrix $X_{i}$ has
independent Gaussian entries distributed according to $\mathcal{N}(
0,\sigma ^{2}/N) $. In this case, $\log \| X_{1}v\| ^{2}$
is distributed in the same way as the random variable
\[
y=\log \Biggl( \frac{\sigma ^{2}}{N}\sum_{k=1}^{N}Y_{k}^{2}\Biggr),
\]%
where $Y_{k}$ are independent standard Gaussian variables. In order to prove
Proposition \ref{Proposition_LD_inequality_unbounded_rv} for this case, it
is enough to show that the following result holds.

\begin{proposition}
Let $y_{i}$ be independent copies of the variable
\[
y=\log \Biggl[ \frac{1}{N}\sum_{k=1}^{N}Y_{k}^{2}\Biggr] .
\]
If $t\leq 1,$ then there exists a function $N_{0}( t) $ such that
for all $N\geq N_{0}( t) $ and all $n,$ the following inequality
holds:%
\[
\Pr \Biggl\{ \Biggl| \frac{1}{n}\sum_{i=1}^{n}y_{i} \Biggr| \geq t\Biggr\} \leq
2e^{-t^{2}nN/8}.
\]
\end{proposition}

\begin{pf} First of all, let us compute
\[
Ee^{yz}=E\Biggl( \frac{1}{N}\sum_{k=1}^{N}Y_{k}^{2}\Biggr) ^{z},
\]%
where $z$ is a real number. By explicit calculation,
\[
E\Biggl( \sum_{i=1}^{N}Y_{i}^{2}\Biggr) ^{z}=\frac{2^{z}\Gamma ( {N%
/2}+z) }{\Gamma ( {N/2}) },
\]%
where $\Gamma ( z) $ is the Gamma function. This formula is valid
for $z>-N/2.$

Let $z=\alpha N,$ where $\alpha >-1/2.$ Then using the Stirling formula for
large $N$, we can write
%
\begin{eqnarray}\label{Stirling_expansion}
E( e^{yz}) &=&E\Biggl( \sum_{i=1}^{N}Y_{i}^{2}\Biggr) ^{z}\nonumber
\\[-8pt]\\[-8pt]
&\sim& \frac{%
1}{\sqrt{1+2\alpha }}N^{\alpha N}\exp \biggl\{ \biggl[ \biggl( \frac{1}{2}%
+\alpha \biggr) \log ( 1+2\alpha ) -\alpha \biggr] N\biggr\} .\nonumber
\end{eqnarray}
Note that for all $\alpha \geq 0,$ $( \frac{1}{2}+\alpha ) \log
( 1+2\alpha ) -\alpha \leq \alpha ^{2},$ and for all $\alpha
>-1/2,$ $( \frac{1}{2}+\alpha ) \log ( 1+2\alpha )
-\alpha \leq 2\alpha ^{2}$ with equalities only for $\alpha =0.$

If $t>0,$ we set $\alpha =t/2$ and $z=( t/2) N,$ and use the fact
that for all sufficiently large $N,$ the asymptotic term in (\ref%
{Stirling_expansion}) dominates all other terms. Hence, we obtain the
estimate
\[
e^{-tz}Ee^{yz}\leq \exp ( -t^{2}N/4) .
\]
If $t\in ( -1,0) $, then we can take $\alpha =t/4$ and $z=(
t/4) N,$ and we obtain
\[
e^{-tz}Ee^{yz}\leq \exp ( -t^{2}N/8) .
\]

By standard arguments, we can translate these inequalities into statements
about probabilities of large deviations. If $0\leq t<1,$ then%
\[
\Pr \Biggl\{ \Biggl| \frac{1}{n}\sum_{i=1}^{n}y_{i}\Biggr| \geq t\Biggr\} \leq
2e^{-t^{2}Nn/8}.
\]
\upqed\end{pf}

\subsection{Matrices with uniformly bounded singular values}

In this section, we are going to prove the second part of Proposition \ref%
{Proposition_LD_inequality_unbounded_rv}. Since $X_{i}$ are i.i.d. and
rotationally invariant, therefore the distribution of $y_{i}=\log (
\| \Pi _{i}v\| ^{2}/\| \Pi _{i-1}v\| ^{2}) $
coincides with the distribution of $\log ( \| X_{1}v\|
^{2}) $ and equals the distribution of the random variable $%
y=\sum_{k=1}^{N}s_{k}u_{k}^{2}.$ Here $u_{k}$ are components of the random
vector $u,$ which is uniformly distributed on the unit sphere and which is
independent of $s_{k}.$\vspace*{1pt}

Let us start with considering large deviations of $x=%
\sum_{k=1}^{N}s_{k}u_{k}^{2}$. Let $\overline{s}^{( N)
}=:N^{-1}\sum_{k=1}^{N}s_{k}.$

\begin{proposition}
Suppose that with probability 1, $| s_{k}| \leq B$ for all $k.$
Then for all $t>0,$
\begin{equation}\label{inequality_probability_of_deviation0}
\quad \max \Biggl[ \Pr \Biggl\{ \sum_{k=1}^{N}s_{k}u_{k}^{2}\leq \overline{s}%
-t\Biggr\} ,\Pr \Biggl\{ \sum_{k=1}^{N}s_{k}u_{k}^{2}\geq \overline{s}%
+t\Biggr\} \Biggr] \leq \exp \biggl\{ -\frac{Nt^{2}}{4B( B+t) }%
\biggr\} .
\end{equation}
\end{proposition}

\begin{pf} Let $x$ denote $\sum_{k=1}^{N}s_{k}u_{k}^{2}$ and let us
estimate $\Pr \{ x\geq \overline{s}^{( N) }+t\} .$ We
will estimate the conditional probability $\Pr \{ x\geq \overline{s}%
^{(N)} +t \vert s_1,\ldots,s_N\} ,$ which we denote as $\Pr \{ x\geq
\overline{s}+t\} $ for simplicity. Note that
\begin{eqnarray*}
\Pr \{ x\geq \overline{s}+t\} &\leq &e^{-z( \overline{s}%
+t) }Ee^{zx} \\
&=&e^{-z( \overline{s}+t) }\bigl( 1+M_{1}z+\tfrac{1}{2!}%
M_{2}z^{2}+\cdots\bigr) ,
\end{eqnarray*}
where $z>0$ and $M_{p}=Ex^{p}.$

Let us use von Neumann's formulas from (\cite{vonneumann41}, pages 373--375)
for the uncentered moments of the random variable $x$. Namely, let
\[
\alpha _{l}=\frac{1}{2l}\sum_{i=1}^{N}( s_{i}) ^{l},
\]%
and let
\[
1+\beta _{1}z+\beta _{2}z^{2}+\beta _{3}z^{3}+\cdots =e^{\alpha _{1}z+\alpha
_{2}z^{2}+\alpha _{3}z^{3}+\cdots }.
\]%
Then, von Neumann's result is that
\[
Ex^{k}=\frac{2^{k}k!}{N( N+2) \cdots ( N+2k-2) }\beta
_{k}.
\]

Using this result, we write
\begin{eqnarray*}
1+M_{1}z+\frac{1}{2!}M_{2}z^{2}+\cdots &=&1+\frac{2}{N}\beta _{1}z+\frac{2^{2}}{%
N( N+2) }\beta _{2}z^{2}+\cdots \\
&\leq &1+\beta _{1}\biggl( \frac{2z}{N}\biggr) +\beta _{2}\biggl( \frac{2z}{N}%
\biggr) ^{2}+\cdots \\
&=&\exp \biggl\{ \alpha _{1}\biggl( \frac{2z}{N}\biggr) +\alpha _{2}\biggl(
\frac{2z}{N}\biggr) ^{2}+\cdots\biggr\} .
\end{eqnarray*}

Next, note that $2\alpha _{1}/N=\overline{s},$ and that $\alpha _{k}\leq
k^{-1}( N/2) B^{k}.$ This implies that
\begin{eqnarray*}
&&e^{-z( \overline{s}+t) }\biggl( 1+M_{1}z+\frac{1}{2!}%
M_{2}z^{2}+\cdots\biggr)
\\
&&\qquad\leq e^{-zt}\exp \biggl\{ \frac{N}{4}\biggl[ \biggl(
\frac{2Bz}{N}\biggr) ^{2}+\biggl( \frac{2Bz}{N}\biggr) ^{3}+\cdots\biggr]
\biggr\} \\
&&\qquad=e^{-zt}\exp \biggl\{ \frac{N}{4}\frac{( {2Bz/N}) ^{2}}{1-%
{2Bz/N}}\biggr\} \\
&&\qquad=e^{-zt}\exp \biggl\{ \frac{B^{2}z^{2}}{N-2Bz}\biggr\} .
\end{eqnarray*}
Let
\[
z_{0}=\frac{Nt}{2B( B+t) }.
\]%
Then
\[
\frac{B^{2}z_{0}^{2}}{N-2Bz_{0}}-z_{0}t=-\frac{Nt^{2}}{4B( B+t) }%
.
\]
Altogether, we get
\[
\Pr \{ x\geq \overline{s}+t\} \leq \exp \biggl\{ -\frac{Nt^{2}}{%
4B( B+t) }\biggr\} .
\]%
The proof of the inequality for $\Pr \{ x\leq \overline{s}-t\} $
is similar.
\end{pf}

\begin{corollary}
Suppose that with probability 1, $s_{k}\leq b\overline{s}$ for all $k.$ Then
for all $t>0,$
\begin{eqnarray}
&&\max \Biggl[ \Pr \Biggl\{ \sum_{k=1}^{N}s_{k}u_{k}^{2}\leq \overline{s}(
1-t) \Biggr\} ,\Pr \Biggl\{ \sum_{k=1}^{N}s_{k}u_{k}^{2}\geq \overline{s%
}( 1+t) \Biggr\} \Biggr]\nonumber
\\
&&\qquad\leq \exp \biggl\{ -\frac{Nt^{2}}{%
4b( b+t) }\biggr\} .
\end{eqnarray}
\end{corollary}

\begin{corollary}
Let $0\leq s_{k}\leq b$ for each $k,$ and $t\in ( 0,1/2) .$ Then:
\newline
\textup{(i)}
%
\begin{equation}\label{inequality_probability_of_deviation1}
\Pr \Biggl\{ \log \sum_{k=1}^{N}s_{k}u_{k}^{2}\geq \log \overline{s}%
+t\Biggr\} \leq \exp \biggl\{ -\frac{Nt^{2}}{4b( b+t) }\biggr\} ;
\end{equation}
\textup{(ii)}
\begin{equation}
\Pr \Biggl\{ \log \sum_{k=1}^{N}s_{k}u_{k}^{2}\geq \log \overline{s}-(
2\log 2) t\Biggr\} \leq \exp \biggl\{ -\frac{Nt^{2}}{4b(
b+t) }\biggr\} ;
\end{equation}
\textup{(iii)}
%
\begin{equation}
\Pr \Biggl\{ \Biggl| \log \sum_{k=1}^{N}s_{k}u_{k}^{2}-\log \overline{s}%
\Biggr| \geq t\Biggr\} \leq 2\exp \biggl\{ -\frac{Nt^{2}}{4c( c+t)
}\biggr\} ,\vspace*{2pt}
\end{equation}
where $c=( 2\log 2) b.$
\end{corollary}

\begin{pf} Let $x$ denote $\sum_{k=1}^{N}s_{k}u_{k}^{2}.$ Then\vspace*{2pt}
\begin{eqnarray*}
\Pr \{ x\geq \overline{s}+t\} &=&\Pr \biggl\{ \log x\geq \log
\overline{s}+\log \biggl( 1+\frac{t}{\overline{s}}\biggr) \biggr\}
\\[2pt]
&\geq &\Pr \biggl\{ \log x\geq \log \overline{s}+\frac{t}{\overline{s}}%
\biggr\} .
\end{eqnarray*}
This and (\ref{inequality_probability_of_deviation0}) prove the first
inequality.  The second inequality is proved similarly, and the third one is
a consequence of the first two inequalities.
\end{pf}

\begin{lemma}
\label{lemma_bound_Laplace_transform1}Suppose that $X$ is a random variable
such that\vspace*{2pt}
\[
\Pr \{ | X| \geq t\} \leq 2\exp \biggl\{ -\frac{Nt^{2}}{%
4c( c+t) }\biggr\} ,
\]%
where $c>0$. Let $| z| <N/(16c).$ Then\vspace*{2pt}
\[
Ee^{zX}\leq \sqrt{32\pi }\sqrt{\frac{c^{2}z^{2}}{N}}\exp \biggl( 2\frac{%
c^{2}z^{2}}{N}\biggr) +3e^{| z| /\sqrt{N}}+2\exp \biggl( -\frac{N}{%
16}\biggr) .
\]
\end{lemma}

\begin{pf} Consider the case when $z\geq 0.$ First, let us estimate $%
\int_{1/\sqrt{N}}^{\infty }e^{zt}\mu ( dt) ,$ where $\mu $ is the
distribution measure of $X.$ Let $F( t) =:\Pr \{ X\geq
t\} .$ Then, by integrating by parts and using the inequalities\vspace*{2pt}
\[
F( t) \leq 2\exp \biggl\{ -\frac{Nt^{2}}{4c( c+t) }%
\biggr\}\vspace*{2pt}
\]
and $N\geq 1,$ we get\vspace*{2pt}
\begin{eqnarray*}
\int_{1/\sqrt{N}}^{\infty }e^{zt}\mu ( dt) &=&F\biggl( \frac{1}{%
\sqrt{N}}\biggr) e^{z/\sqrt{N}}+z\int_{1/\sqrt{N}}^{\infty }e^{zt}F(
t)\, dt
\\[2pt]
&\leq &2e^{-1/[ 4c( c+1) ] }e^{z/\sqrt{N}}+2z\int_{1/%
\sqrt{N}}^{\infty }e^{zt}\exp \biggl\{ -\frac{Nt^{2}}{4c( c+t) }%
\biggr\} \, dt.
\end{eqnarray*}

\noindent
In order to estimate the integral in the last line, we divide it into two
pieces, $\int_{1/\sqrt{N}}^{b}$ and $\int_{b}^{\infty }.$ Then
\begin{eqnarray*}
\int_{1/\sqrt{N}}^{b}e^{zt}\exp \biggl\{ -\frac{Nt^{2}}{4c( c+t) }%
\biggr\} \,dt &\leq &\int_{-\infty }^{\infty }e^{zt}\exp \biggl\{ -\frac{Nt^{2}%
}{8c^{2}}\biggr\}\, dt
 \\
&=&\exp \biggl( \frac{2c^{2}}{N}z^{2}\biggr) \int_{-\infty }^{\infty }\exp
\biggl\{ -\frac{N}{8c^{2}}\biggl( t-\frac{4c^{2}}{N}z\biggr) ^{2}\biggr\}\, dt
\\
&=&\sqrt{\frac{8\pi c^{2}}{N}}\exp \biggl( \frac{2c^{2}}{N}z^{2}\biggr) .
\end{eqnarray*}
Next, for the second piece, we have
\begin{eqnarray*}
\int_{b}^{\infty }e^{zt}\exp \biggl\{ -\frac{Nt^{2}}{4c( c+t) }%
\biggr\} \, dt &\leq &\int_{b}^{\infty }e^{zt}\exp \biggl\{ -\frac{Nt}{8c}%
\biggr\} \,dt
\\
&=&\frac{1}{N/( 8c) -z}\exp \biggl( -\biggl( \frac{N}{8c}-z\biggr)
c\biggr) \\
&\leq &\frac{16c}{N}\exp \biggl( -\frac{N}{16}\biggr) ,
\end{eqnarray*}
where we used the assumption that $z\leq N/( 16c) .$

Hence, combining the previous inequalities and using the assumption that $%
z\leq N/( 16c) $ again, we get
\begin{eqnarray*}
\int_{1/\sqrt{N}}^{\infty }e^{zt}\mu ( dt) &\leq &2e^{-1/[
4c( c+1) ] }e^{z/\sqrt{N}}+2z\sqrt{\frac{8\pi c^{2}}{N}}%
\exp \biggl( \frac{2c^{2}}{N}z^{2}\biggr)
\\
&&{}+2\exp \biggl( -\frac{N}{16}\biggr) .
\end{eqnarray*}

In addition,
\[
\int_{-\infty }^{1/\sqrt{N}}e^{zt}\mu ( dt) \leq e^{z/\sqrt{N}}.
\]

Combining all the parts, we get
\begin{eqnarray*}
\int_{-\infty }^{\infty }e^{zt}\mu ( dt) &\leq &\sqrt{\frac{32\pi
c^{2}z^{2}}{N}}\exp \biggl( \frac{2c^{2}}{N}z^{2}\biggr)
+\bigl( 1+2e^{-1/[ 4c( c+1) ] }\bigr) e^{z/\sqrt{N}}
\\
&&{}+2\exp \biggl( -\frac{N}{16}\biggr) ,
\end{eqnarray*}
from which the claim of the lemma follows for $z\geq 0$. The case when $%
z\leq 0$ is similar.
\end{pf}

\begin{corollary}
Let $X=\log ( \sum_{k=1}^{N}s_{k}u_{k}^{2}) -\log (
\overline{s}) $ and let $| z| \leq N/(16c),$ where $%
c=( 2\log 2) b.$ Then
\[
Ee^{zX}\leq \sqrt{32\pi }\sqrt{\frac{c^{2}z^{2}}{N}}\exp \biggl( 2\frac{%
c^{2}z^{2}}{N}\biggr) +3e^{| z| /\sqrt{N}}+2\exp \biggl( -\frac{N}{%
16}\biggr) .
\]
\end{corollary}

\begin{pf} This follows directly from Lemma \ref%
{lemma_bound_Laplace_transform1} and inequality (\ref%
{inequality_probability_of_deviation1}).
\end{pf}

\begin{pf*}{Proof of the second part of Proposition \protect\ref%
{Proposition_LD_inequality_unbounded_rv}} Note that
\[
\log \| \Pi _{n}u\| ^{2}=\sum_{i=1}^{n}\log \Biggl[
\sum_{k=1}^{N}s_{k}^{( i,N) }\bigl( u_{k}^{( i,N)
}\bigr) ^{2}\Biggr],
\]
where $u_{k}^{( i,N) }$ are components of independent
Haar-distributed $N$-vectors $u^{( i,N) }.$ Let
\[
Y_{i}=\log \Biggl( \sum_{k=1}^{N}s_{k}^{( i,N) }\bigl(
u_{k}^{( i,N) }\bigr) ^{2}\Biggr) -\log \bigl( \overline{s}%
^{( i,N) }\bigr) .
\]%
We aim to estimate
\[
\Pr \Biggl\{ \Biggl| \sum_{i=1}^{n}Y_{i}\Biggr| >nt\Biggr\} .
\]
As usual,
\[
\Pr \Biggl\{ \sum_{i=1}^{n}Y_{i}>nt\Biggr\} \leq e^{-nzt}(
Ee^{zY_{i}}) ^{n},
\]%
where $z>0.$

Note that by Assumption \ref{AB},
\[
s_{k}^{( i,N) }\leq b\overline{s}^{( i,N) },
\]
hence, our previous lemmas are applicable.

We set $z=tN/( 4c^{2}) $ and assume that $N\geq 4/t^{2}.$ (Note
that the assumption that $t\in ( 0,1/4] $ implies that $z\leq
N/( 16c) .$) Then, by the previous lemma, we have
\begin{eqnarray*}
Ee^{zY_{i}} &\leq &\sqrt{32\pi }\sqrt{\frac{c^{2}z^{2}}{N}}\exp \biggl( 2%
\frac{c^{2}z^{2}}{N}\biggr) +3e^{z/\sqrt{N}}+2\exp \biggl( -\frac{N}{16}\biggr)
\\
&\leq &\sqrt{2\pi }\sqrt{\frac{t^{2}N}{c^{2}}}\exp \biggl( \frac{1}{8}\frac{%
t^{2}N}{c^{2}}\biggr) +3\exp \biggl( \frac{t\sqrt{N}}{4c^{2}}\biggr) +2\exp
\biggl( -\frac{N}{16}\biggr).
\end{eqnarray*}
Since $N\geq 4/t^{2},$ then the first term dominates the other two terms,
and we can write
\[
Ee^{zY_{i}}\leq \Biggl( \sqrt{2\pi }\sqrt{\frac{t^{2}N}{c^{2}}}+5\Biggr) \exp
\biggl( \frac{t^{2}N}{8c^{2}}\biggr) .
\]
Hence,
\begin{eqnarray*}
e^{-nzt}( Ee^{zY_{i}}) ^{n} &\leq &\Biggl( \sqrt{2\pi }\sqrt{\frac{%
t^{2}N}{c^{2}}}+5\Biggr) ^{n}\exp \biggl( -\frac{t^{2}N}{8c^{2}}n\biggr)
\\
&=&\exp \bigl\{ -n\bigl[ -\log \bigl( \sqrt{2\pi }( t/c) \sqrt{N}%
+5\bigr) +( t^{2}/8c^{2}) N\bigr] \bigr\} .
\end{eqnarray*}
Clearly, we can find an $N_{0}( t) $ such that for all $%
N>N_{0}( t) $
\[
e^{-nzt}( Ee^{zY_{i}}) ^{n}\leq \exp \{ -n(
t^{2}/16c^{2}) N\} .
\]%
Hence, for all $N>N_{0}( t) $
\begin{equation}\label{inequality_probability_of_deviation2}
\Pr \Biggl\{ \frac{1}{n}\sum_{i=1}^{n}\bigl[ y_{i}^{( N) }-\log
\bigl( \overline{s}^{( i,N) }\bigr) \bigr] >t\Biggr\} \leq \exp
\{ -n( t^{2}/16c^{2}) N\} .
\end{equation}
The case of the inequality
\begin{equation}
\Pr \Biggl\{ \frac{1}{n}\sum_{i=1}^{n}\bigl[ y_{i}^{( N) }-\log
\bigl( \overline{s}^{( i,N) }\bigr) \bigr] <-t\Biggr\} \leq
\exp \{ -n( t^{2}/16c^{2}) N\}
\end{equation}
is similar. Finally, note that $16c^{2}\leq 32b^{2}.$
\end{pf*}

\section{Necessary condition}\label{section_necessity}

Let us introduce the following assumption.

\renewcommand{\theassumption}{\Alph{assumption}$'$}
\setcounter{assumption}{3}
\begin{assumption}\label{AD'}
 $E[ \log \| X_{i}^{( N) }u\| %
] ^{2}$ exists and bounded by a constant that does not depend on $N.$
\end{assumption}

\begin{theorem}
Let Assumptions \ref{AA}, \ref{AB} and \ref{AD'} hold. Suppose that for every $%
\delta >0$ there exists such an $n_{0}( \delta ) $ that
\begin{equation}
\Pr \bigl\{ \bigl| n^{-1}\log \| \Pi _{n}\| -E\log \bigl\|
X_{1}^{( N) }v\bigr\| \bigr| \geq \delta \bigr\} \leq \delta
\end{equation}
for all $N$ and all $n\geq n_{0}( \delta ) .$ Let $b(
N) $ is an arbitrary function of $N$ such that $\lim_{N\rightarrow
\infty }b( N) =+\infty .$ Then
\[
\lim_{N\rightarrow \infty }\Pr \bigl\{ \bigl\| X_{1}^{( N)
}\bigr\| \geq b( N) \bigr\} =0.
\]
\end{theorem}

\begin{pf} Let $v_{0}$ be such a unit vector that $\|
X_{1}^{( N) }\| =\| X_{1}^{( N)
}v_{0}\| .$ Note that $X_{1}^{( N) }v_{0}$ has the Haar
distribution by assumption of rotational invariance. By using the fact that $%
\| \Pi _{n}\| \geq \| \Pi _{n}v_{0}\| ,$ we can write
the inequality
\[
n^{-1}\log \| \Pi _{n}\| \geq \frac{\log \| X_{1}^{(
N) }\| }{n}+\frac{1}{n}\sum_{i=2}^{n}\log \bigl( X_{i}^{(
N) }u_{i}\bigr) ,
\]%
where $u_{i}$ are independent Haar-distributed vectors. By using Assumption~\ref{AD'}, we~can conclude that $n^{-1}\sum_{i=2}^{n}\log (
X_{i}^{( N) }u_{i}) $ converges in probability to\break $E \log
\| X_{1}^{( N) } v\| $ and that the convergence is
uniform in $N.$ This fact and the supposition of the theorem imply that $%
n^{-1}\log \| X_{1}^{( N) }\| $ must converge in
probability to zero as $n\rightarrow \infty ,$ and that the convergence must
be uniform in $N.$ If the conclusion of the theorem were invalid, then for
some $\delta >0$ and all $n,$ we could find an $N=N( n,\delta ) $
such that $\Pr \{ \log \| X_{1}^{( N) }\| \geq
n\delta \} \geq \delta ,$ and this would contradict the uniform
convergence of $n^{-1}\log \| X_{1}^{( N) }\| $ to
zero.
\end{pf}

\section{Conclusion}\label{section_conclusion}

In this paper, we found sufficient conditions that ensure that the
convergence rate in the Furstenberg--Kesten theorem is uniform with respect
to the dimension of the space in which matrices operate. Let us call this
phenomenon dimensional uniformity of convergence.

Several interesting questions remain to be answered. First, is it possible
to prove the dimensional uniformity of convergence for random matrices which
are not rotationally invariant, for example, for Wigner matrices?

Second, assuming rotational invariance, what characterizes the laws of
singular values $s_{k}^{( i,N) },$ for which the dimensional
uniformity of convergence holds? In other words, what are necessary and
sufficient conditions for dimensional uniformity of convergence?


\printaddresses

\end{document}